\documentclass[12pt]{article}
\usepackage{amssymb,epsfig,psfig}

\begin{document}

\title{Exact Polynomial Eigenmodes for Homogeneous Spherical 3-Manifolds}

\author{Jeff Weeks}

\maketitle

\begin{abstract}
    \noindent Observational data hint at a finite universe,
    with spherical manifolds such as the Poincar\'e
    dodecahedral space tentatively providing the best fit.
    Simulating the physics of a model universe
    requires knowing the eigenmodes of the Laplace operator
    on the space.  The present article provides explicit polynomial
    eigenmodes for all globally homogeneous 3-manifolds:
    the Poincar\'e dodecahedral space $S^3/I^*$, the binary
    octahedral space $S^3/O^*$, the binary tetrahedral space
    $S^3/T^*$, the prism manifolds $S^3/D_m^*$
    and the lens spaces $L(p,1)$.
\end{abstract}

\section{Introduction}
\label{SectionIntroduction}

The past decade has seen intense work on multiconnected
3-manifolds as models for the physical universe. Well-proportioned
3-manifolds explain the missing broad-scale fluctuations in the
cosmic microwave background, first discovered by the COBE
satellite~\cite{hin96} and later confirmed by the WMAP
satellite~\cite{ben03,spergel}.  (Please see Ref.~\cite{wee04} for
an elementary exposition.)  A {\it well-proportioned} 3-manifold
is one whose three dimensions are of comparable magnitudes.
Ill-proportioned manifolds, with one dimension significantly
larger or smaller than the other two, fail to explain the missing
broad-scale fluctuations and indeed predict exactly the opposite,
namely elevated broad-scale fluctuations~\cite{WellProportioned},
contrary to observations. Current work, therefore, focuses on
well-proportioned spaces.

The density of ordinary matter alone would suggest a hyperbolic
universe, and so ten years ago researchers studied hyperbolic
models. The situation changed dramatically in late 1998 with the
discovery of a still-mysterious vacuum energy that raises the
universe's mass-energy density parameter~$\Omega$ to the level
required for a flat space ($\Omega = 1$) or a slightly spherical
space ($\Omega > 1$).  The first-year WMAP paper \cite{ben03} put
$\Omega$ at $1.02 \pm 0.02$ at the $1\sigma$~level, while the
three-year WMAP paper \cite{spergel} reports six different
distributions ranging from a rather flat $\Omega =
1.003^{+0.017}_{-0.013}$ to a surprisingly curved $\Omega =
1.037^{+0.015}_{-0.021}$ depending on what external data set one
uses to resolve the geometrical degeneracy.

The Poincar\'e dodecahedral space \cite{WeberSeifert33}, defined
as the quotient $S^3/I^*$ of the 3-sphere by the binary
icosahedral group $I^*$, explains both the missing fluctuations
and the observed mass-energy density~\cite{Nature} and so
researchers are now modelling it more precisely for better
comparison to observations. State-of-the-art simulations find that
the Poincar\'e dodecahedral space matches observed broad-scale
fluctuations when $1.015 < \Omega < 1.020$~\cite{aurich,aurich2}
or $1.022 < \Omega < 1.034$~\cite{gundermann}, in excellent
agreement with the three-year WMAP estimates of $\Omega$. However,
other topologies, notably the quotient $S^3/O^*$ of the 3-sphere
modulo the binary octahedral group, also remain
viable~\cite{aurich2,gundermann}.

Attempts to confirm a multiconnected universe using the
circles-in-the-sky method~\cite{cornish} have
failed~\cite{cornish2}.  (Exception:  one research group claims to
have found hints of matching circles~\cite{roukema} and a second
group has confirmed the match~\cite{key}, but the former provide
no statistical analysis while the latter find the match to lie
below the false positive threshold, so the former's claim remains
unconvincing.)  Accepting the result that any potential circle
pairs are undetectable~\cite{cornish2}, the question remains:  are
the circles really not there, or are they merely hidden by various
sources of contamination, such as the Doppler and integrated
Sachs-Wolfe components of the microwave background?  Answering
this question requires great care, because the level of
contamination depends both on the topology and on the choice of
cosmological parameters. Results so far remain inconsistent:  the
circle-searchers' own analysis finds their negative results to be
robust for a dodecahedral universe in spite of the
contamination~\cite{key}, while other researchers find
contamination strong enough to hide matching
circles~\cite{aurich3}.
\\

To determine the observational consequences of a given
cosmological model, researchers take a Fourier approach and
express physical quantities, such as density fluctuations in the
primordial plasma, as linear combinations of the eigenmodes of the
Laplacian (more briefly, the {\it modes}), which can then be
integrated forward in time.  Thus all studies of cosmic topology
require knowing the modes of the 3-manifold under consideration.
Different research groups have determined the modes in different
ways~\cite{aurich,aurich2,gundermann,lachiezerey,lwugl}, but so
far all approaches have required some sort of numerical
computation, either the extraction of the eigenvectors of a large
matrix or the solution of a large set of simultaneous equations.

The present article provides the modes as explicit polynomials
with integer coefficients. That is, for each manifold $S^3/\Gamma$
and each wavenumber $k$, we provide a finite set of
$\Gamma$-invariant polynomials of degree $k$ spanning the full
space of modes.  We provide these polynomials for the binary
icosahedral space $S^3/I^*$ (better known as the Poincar\'e
dodecahedral space), the binary octahedral space $S^3/O^*$, the
binary tetrahedral space $S^3/T^*$, the binary dihedral spaces
$S^3/D_m^*$ (better known as prism manifolds), the homogeneous
lens spaces $L(p,1)$, and the 3-sphere itself.  These spaces
comprise the full set of globally homogeneous spherical
3-manifolds (called {\it single action spaces} in the
classification of Ref.~\cite{glluw}).

The ideas in the present article draw heavily on Klein's 1884 {\it
Vorlesungen \"uber das Ikosaeder} \cite{Klein}, extending Klein's
work to produce full bases of explicit polynomials.

Sections~\ref{SectionLiftingS2toS3} through
\ref{SectionHarmonicity} are elementary in nature, laying a
foundation, establishing terminology, and translating into
geometric language some concepts that have recently appeared in
the cosmic topology literature only in quantum mechanical bracket
language.  Sections~\ref{SectionModesS3I} through
\ref{SectionModesS3} provide the real content of this article,
namely the explicit polynomials for the modes of $S^3/I^*$,
$S^3/O^*$, $S^3/T^*$, $S^3/D_m^*$, $L(p,1)$ and $S^3$.
Section~\ref{SectionFutureWork} outlines future work.

\section{Lifting $S^2$ to $S^3$}
\label{SectionLiftingS2toS3}

Let us parameterize the 2-sphere $S^2$ first as the Riemann sphere
$\mathbb C \cup \{\infty\}$ and then, more elegantly, as the
complex projective line $\mathbb C P^1$.  First map the unit
2-sphere onto the equatorial plane via stereographic projection
from the south pole:  typical points map via $(x,y,z) \mapsto
\frac{x}{z+1} + i\frac{y}{z+1}$ while the south pole $(0,0,-1)$
maps to $\infty$. To accommodate the south pole into the same
format as the typical points, write each image point as a {\it
formal fraction}.  Each typical point maps to a formal fraction
$\frac{x + i y}{z+1}$, while the south pole maps to the formal
fraction $\frac{1}{0}$. Two formal fractions
$\frac{\alpha}{\beta}$ and $\frac{\alpha'}{\beta'}$ are equivalent
if and only if $\alpha \beta' = \alpha' \beta$.  Let
$[\frac{\alpha}{\beta}]$ denote the equivalence class of all
formal fractions $\frac{\alpha'}{\beta'}$ equivalent to
$\frac{\alpha}{\beta}$.  The set of all such equivalence classes
defines the {\it complex projective line} $\mathbb C P^1$. The
indeterminate fraction $\frac{0}{0}$ is of course excluded from
the discussion. In summary, the set $\mathbb C P^1$ of equivalence
classes $[\frac{\alpha}{\beta}]$ of of nontrivial complex formal
fractions parameterizes the 2-sphere $S^2$.

Writing the formal fraction $\frac{\alpha}{\beta}$ as an ordered
pair $(\alpha,\beta) \in \mathbb C^2$ immediately yields a map
$\mathbb C^2 - (0,0) \rightarrow \mathbb C P^1$ sending each point
$(\alpha,\beta) \in \mathbb C^2$ to the corresponding equivalence
class $[\frac{\alpha}{\beta}] \in \mathbb C P^1$. The radial
direction in $\mathbb C^2$ is largely irrelevant, because
$(\alpha,\beta)$ and $(r \alpha, r\beta)$ map to the same class
$[\frac{\alpha}{\beta}] = [\frac{r \alpha}{r \beta}] \in \mathbb C
P^1$ for all real $r \neq 0$. We may therefore restrict the map's
domain to the unit 3-sphere $S^3 \subset \mathbb C^2$ without
compromising its image in $\mathbb C P^1$. Recalling from the
preceding paragraph that $\mathbb C P^1$ parameterizes $S^2$, this
yields a map $S^3 \rightarrow S^2$.

The {\it Clifford flow}\, on $S^3 \subset \mathbb C^2$ is the
fixed-point-free motion taking $(\alpha,\beta) \mapsto (e^{i t}
\alpha, e^{i t} \beta)$ for $0 \leq t \leq 2\pi$.  Under the
Clifford flow, each point $(\alpha,\beta)$ traces out a great
circle, which collectively form a set of {\it Clifford parallels}.
The points along a given Clifford parallel all map to the same
equivalence class $[\frac{e^{i t} \alpha}{e^{i t} \beta}] =
[\frac{\alpha}{\beta}] \in \mathbb C P^1$.  In other words, the
preimage of each point of $S^2$ is a Clifford parallel in $S^3$
which we will call a {\it fiber} of the map.

\section{Constructing the Groups}
\label{SectionGroupConstruction}

\begin{figure}
\centerline{\psfig{file=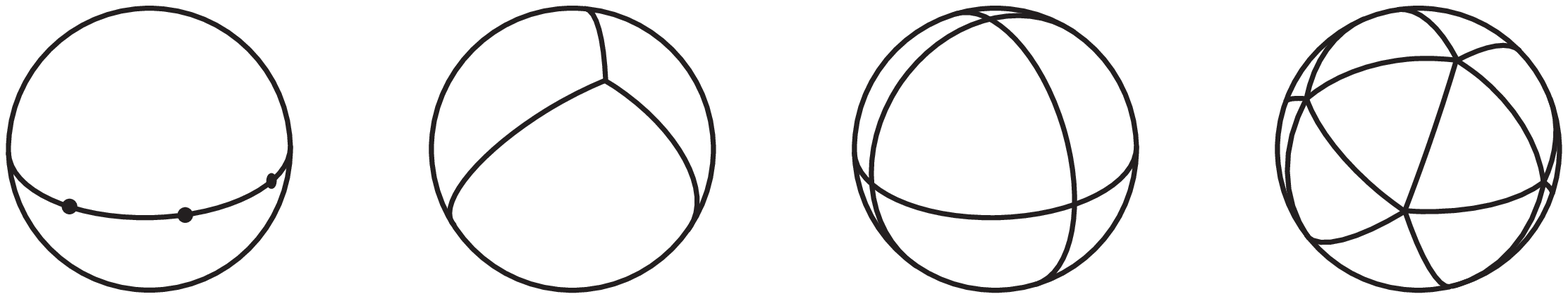, width=12cm}} \caption{The
dihedron, tetrahedron, octahedron and icosahedron project radially
to tilings of the 2-sphere.  The dihedron has only two faces, one
filling the whole northern hemisphere and the other filling the
whole southern hemisphere.  The dihedron's vertices and edges lie
evenly along the equator.  (Artwork by Adam Weeks Marano.)}
\label{FigureRegularPolyhedra}
\end{figure}

The classic platonic solids project radially to tilings of $S^2$
(Figure~\ref{FigureRegularPolyhedra}).  We ignore the cube and
dodecahedron, because they are dual to the octahedron and
icosahedron respectively, but include the dihedron, which works
out fine as a tiling of $S^2$ even though it falls flat as a
traditional polyhedron.  The orientation-preserving symmetries of
these tilings comprise the {\it dihedral group} $D_m$, the {\it
tetrahedral group} $T$, the {\it octahedral group} $O$ and the
{\it icosahedral group} $I$, of order $2m$, 12, 24 and 60,
respectively.

Every matrix $( \phantom{|}^a_c \phantom{|}^b_d \phantom{|} ) \in
GL_2 \mathbb C$ respects Clifford parallels in the sense that
left-multiplication takes each fiber $( e^{i t} \alpha, e^{i t}
\beta )$ to another fiber $( e^{i t} (a \alpha + b \beta), e^{i t}
(c \alpha + d \beta) )$. Therefore the matrix action projects down
to a well-defined map $S^2 \rightarrow S^2$ taking
$[\frac{\alpha}{\beta}]$ to $[\frac{a \alpha + b \beta}{c \alpha +
d \beta}]$.  Restricting our attention to unitary matrices $(
\phantom{|}^a_c \phantom{|}^b_d \phantom{|} ) \in U(2)$ ensures
rigid motions of $\mathbb C^2$, which in turn project down to
rigid motions of $S^2$.  The unitary matrices, however, still
allow room for ``sliding along the fibers'', for example via the
matrix $( \phantom{|}^{e^{i\theta}}_0 \phantom{|}^0_{e^{i\theta}}
\phantom{|} )$, which of course has no effect when projected down
to $S^2$.  To obtain an (almost) unique matrix for each rotation
of $S^2$, restrict to the special unitary group $SU(2)$, whose
matrices take the form $( \phantom{|}^z_w \phantom{|}^{-\bar
w}_{\phantom{-}\bar z} \phantom{|} )$.  Each rotation of $S^2$ is
realized by exactly two special unitary matrices $\pm(
\phantom{|}^z_w \phantom{|}^{-\bar w}_{\phantom{-}\bar z}
\phantom{|} )$.

How may we construct the pair of matrices $\pm( \phantom{|}^z_w
\phantom{|}^{-\bar w}_{\phantom{-}\bar z} \phantom{|} )$ realizing
a given rotation of $S^2$?  To construct an order $n$ rotation
about a desired fixed point $[\frac{p}{q}]$, we require that
\begin{equation}
  \left(
    \begin{array}{cc}
      z & -\bar w \\
      w & \phantom{-}\bar z \\
    \end{array}
  \right)
  \left(
    \begin{array}{c}
      p \\
      q \\
    \end{array}
  \right)
  =
  \left(
    \begin{array}{c}
      e^{\frac{i\pi}{n}} p \\
      e^{\frac{i\pi}{n}} q \\
    \end{array}
  \right).
\end{equation}
The phase factor $e^{\frac{i\pi}{n}}$ ensures that we rotate the
correct amount about the given fixed point, so the isometry will
have the desired order $n$.  A quick calculation gives
\begin{equation}
  \label{EqnRotationMatrix}
  \left(
    \begin{array}{c}
      z \\
      \bar w \\
    \end{array}
  \right)
  =
  \pm
  \frac{1}{p \bar p + q \bar q}
  \left(
    \begin{array}{cc}
      \phantom{-}\bar p & q \\
               - \bar q & p \\
    \end{array}
  \right)
  \left(
    \begin{array}{c}
      e^{\frac{+i\pi}{n}}      p \\
      e^{\frac{-i\pi}{n}} \bar q \\
    \end{array}
  \right)
\end{equation}
and $\pm( \phantom{|}^z_w \phantom{|}^{-\bar w}_{\phantom{-}\bar
z} \phantom{|} )$ follows immediately.

\begin{table}
\begin{center}
\begin{tabular}{|c|c|c|}
  \hline
  center & order & matrices \\
  \hline
  - & 1 &
  $\pm\left(
    \begin{array}{cc}
      1 & \phantom{-}0 \\
      0 & \phantom{-}1 \\
    \end{array}
  \right)$
  \\
  $0,\infty$ & 2 &
  $\pm\left(
    \begin{array}{cc}
      i &  \phantom{-}0 \\
      0 &      -i       \\
    \end{array}
  \right)$
  \\
  $\pm 1$ & 2 &
  $\pm\left(
    \begin{array}{cc}
      0 & \phantom{-}i \\
      i & \phantom{-}0 \\
    \end{array}
  \right)$
  \\
  $\pm i$ & 2 &
  $\pm\left(
    \begin{array}{cc}
       0 &     -1       \\
       1 & \phantom{-}0 \\
    \end{array}
  \right)$
  \\
  \hline
\end{tabular}
\end{center}
\caption{The four rotations in $D_2 \subset SO(3)$ correspond to
eight matrices in $D_2^* \subset SU(2)$.}
\label{TableElementsD2}
\end{table}

Equation~(\ref{EqnRotationMatrix}) makes it easy to write down
matrices for the groups $D_m$, $T$, $O$ and $I$.  Two matrices
realize each rotation of $S^2$, so for example the 4-element group
$D_2 \subset SO(3)$ is realized by an 8-element matrix group
$D_2^* \subset SU(2)$ called the {\it binary dihedral group of
order 8} (Table~\ref{TableElementsD2}).  One might hope to extract
a 4-element subgroup of $D_2^*$ realizing $D_2$ directly, but this
is impossible because squaring an element of order 2 always gives
$( \phantom{|}^{-1}_{\phantom{-}0} \phantom{|}^{\phantom{-}0}_{-1}
\phantom{|} )$, never $( \phantom{|}^{1}_{0} \phantom{|}^{0}_{1}
\phantom{|} )$, and once $( \phantom{|}^{-1}_{\phantom{-}0}
\phantom{|}^{\phantom{-}0}_{-1} \phantom{|} )$ is in a subgroup,
so is the negative of every matrix in that subgroup.

Still using Equation~(\ref{EqnRotationMatrix}), one may easily
write down matrices for the 12-element tetrahedral group $T
\subset SO(3)$, giving the 24-element {\it binary tetrahedral
group} $T^* \subset SU(2)$, and similarly for the 48-element {\it
binary octahedral group} $O^*$ and the 120-element {\it binary
icosahedral group} $I^*$.  The quotient $S^3/I^*$ defines the
Poincar\'e dodecahedral space.

\section{Constructing Symmetric Polynomials}
\label{SectionSymmetricPolynomials}

Parameterize the 2-sphere $S^2$ as $\mathbb CP^1$
(Section~\ref{SectionLiftingS2toS3}) and let $P =
\{[\frac{p_1}{q_1}],\dots,[\frac{p_n}{q_n}]\}$ be a set of points
thereon, whose symmetry we hope to capture in a polynomial.  As a
starting point, the polynomial
\begin{equation}
  (z - \frac{p_1}{q_1})
  \cdots
  (z - \frac{p_n}{q_n})
\end{equation}
has roots exactly at $P$.  Replacing the variable $z$ with a
formal fraction $[\frac{\alpha}{\beta}] \in \mathbb CP^1$ gives
\begin{equation}
  (\frac{\alpha}{\beta} - \frac{p_1}{q_1})
  \cdots
  (\frac{\alpha}{\beta} - \frac{p_n}{q_n})
\end{equation}
which Klein, in his 1884 {\it Vorlesungen \"uber das Ikosaeder}
\cite{Klein}, writes as a homogeneous polynomial
\begin{equation}
\label{EqnGeneralPolynomial}
  (q_1 \alpha - p_1 \beta)
  \cdots
  (q_n \alpha - p_n \beta).
\end{equation}
I thank Peter Kramer for his recent article~\cite{Kramer} pointing
out the relevance of Klein's work to current investigations in
cosmic topology.

To test whether the polynomial~(\ref{EqnGeneralPolynomial}) is
invariant under a symmetry $\gamma =( \phantom{|}^z_w
\phantom{|}^{-\bar w}_{\phantom{-}\bar z} \phantom{|} ) \in
SU(2)$, rewrite expression~(\ref{EqnGeneralPolynomial}) as
\begin{equation}
  \prod_{i=1}^n \quad
  \left(
    \begin{array}{cc}
      p_i & q_i \\
    \end{array}
  \right)
  \left(
    \begin{array}{rr}
      0 & -1 \\
      1 &  0 \\
    \end{array}
  \right)
  \left(
    \begin{array}{c}
      \alpha \\
      \beta  \\
    \end{array}
  \right)
\end{equation}
and let $\gamma^{-1} = ( \phantom{|}^{\phantom{-}\bar z}_{-w}
\phantom{|}^{\bar w}_{z} \phantom{|} )$ act on $(\alpha, \beta)$,
transforming the polynomial to
\begin{equation}
  \prod_{i=1}^n \quad
  \left(
    \begin{array}{cc}
      p_i & q_i \\
    \end{array}
  \right)
  \left(
    \begin{array}{rr}
      0 & -1 \\
      1 &  0 \\
    \end{array}
  \right)
  \left(
    \begin{array}{rr}
      \bar z & \bar w \\
        -w   &   z    \\
    \end{array}
  \right)
  \left(
    \begin{array}{c}
      \alpha \\
      \beta  \\
    \end{array}
  \right)
\end{equation}
which equals
\begin{equation}
  \prod_{i=1}^n \quad
  \left[
  \left(
    \begin{array}{cc}
      p_i & q_i \\
    \end{array}
  \right)
  \left(
    \begin{array}{rr}
         z    &   w    \\
      -\bar w & \bar z \\
    \end{array}
  \right)
  \right]
  \left(
    \begin{array}{rr}
      0 & -1 \\
      1 &  0 \\
    \end{array}
  \right)
  \left(
    \begin{array}{c}
      \alpha \\
      \beta  \\
    \end{array}
  \right).
\end{equation}
In effect the action of $\gamma$ transforms each root $(p_i, q_i)$
according to the rule $(p_i, q_i) \mapsto (p_i, q_i)(
\phantom{|}^{\phantom{-}z}_{-\bar w} \phantom{|}^w_{\bar z}
\phantom{|} ) = (\gamma)( \phantom{|}^{p_i}_{q_i} \phantom{|} )$.
\\

\noindent {\bf Example \ref{SectionSymmetricPolynomials}.1.}
Consider the points $P = \{\frac{1}{1}, \frac{i}{1},
\frac{-1}{\phantom{-}1}, \frac{-i}{\phantom{-}1}\}$ with
polynomial
\begin{equation}
\label{EqnToyPolynomial}
  (\alpha - \beta)(\alpha - i \beta)(\alpha + \beta)(\alpha + i\beta)
  = \alpha^4 - \beta^4,
\end{equation}
and consider the matrix $\gamma_2 = ( \phantom{|}^{i}_{0}
\phantom{|}^{\phantom{-}0}_{-i} \phantom{|} ) \in SU(2)$.
Geometrically, $\gamma_2$ projects down to an order~2 rotation of
$S^2$ about the north pole, which interchanges the points of $P$
in pairs.  Acting on $S^3 \subset C^2$, $\gamma_2$ takes the roots
to
\begin{equation}
\label{EqnTransformedToyRoots}
  \gamma P
  = \{\frac{\phantom{-}i}{-i}, \frac{-1}{-i},
      \frac{-i}{-i}, \frac{\phantom{-}1}{-i}\}
  = \{-i\frac{-1}{\phantom{-}1}, -i\frac{-i}{\phantom{-}1},
      -i\frac{1}{1}, -i\frac{i}{1}\},
\end{equation}
with formal fractions treated as vectors, so for example
$-i\frac{-1}{\phantom{-}1} = \frac{(-i)(-1)}{(-i)(\phantom{-}1)}$.
The action of $\gamma_2$ permutes the roots and also multiplies
them by a factor of $-i$.  Fortunately these changes leave the
polynomial invariant,
\begin{equation}
\label{EqnTransformedToyPolynomial}
  (-i)(\alpha + \beta)\;(-i)(\alpha + i\beta)\;
  (-i)(\alpha - \beta)\;(-i)(\alpha - i \beta)
  = \alpha^4 - \beta^4.
\end{equation}
\\

\noindent {\bf Example \ref{SectionSymmetricPolynomials}.2.}
Consider the same points $P = \{\frac{1}{1}, \frac{i}{1},
\frac{-1}{\phantom{-}1}, \frac{-i}{\phantom{-}1}\}$ with the same
polynomial $\alpha^4 - \beta^4$ as in
Example~\ref{SectionSymmetricPolynomials}.1, but now act by
$\gamma_4 = ( \phantom{|}^{e^{i\pi/4}}_{\phantom{-}0}
\phantom{|}^{\phantom{-}0}_{e^{-i\pi/4}} \phantom{|} )$.
Geometrically, $\gamma_4$ projects down to an order~4 rotation of
$S^2$ about the north pole, which permutes the points of $P$
cyclically.  Acting on $S^3 \subset C^2$, $\gamma_4$ takes the
roots to
\begin{equation}
\label{EqnTransformedToyRoots2}
  \gamma_4 P
  = \{e^{-i\pi/4}\frac{i}{1}, e^{-i\pi/4}\frac{-1}{\phantom{-}1},
      e^{-i\pi/4}\frac{-i}{\phantom{-}1}, e^{-i\pi/4}\frac{1}{1}\}.
\end{equation}
Again the roots have been permuted, but this time they are
multiplied by a factor of $e^{-i\pi/4}$, so the polynomial maps to
\begin{equation}
\label{EqnTransformedToyPolynomial2}
  (e^{-i\pi/4})^4\;(\alpha - i \beta)(\alpha + \beta)(\alpha + i\beta)(\alpha - \beta)
  = -(\alpha^4 - \beta^4).
\end{equation}
Thus $\gamma_4$ does not leave this polynomial invariant, but
rather sends it to its negative.

\section{Generalized Complex Derivatives}
\label{SectionGeneralizedDerivatives}

Let us extend the complex derivative operator
$\frac{\partial}{\partial z}$ from the class of
complex-differentiable (i.e. analytic) functions to the broader
class of real-differentiable (i.e. smooth) functions.

A complex-valued function $f$ of a complex variable $z = x + i y$
is differentiable in the complex sense if and only if it is
differentiable in the real sense (as a function of $x$ and $y$)
and, in addition,
\begin{equation}
\label{EqnEqualPartials}
  \frac{\partial f}{\partial y}
  =
  i \frac{\partial f}{\partial x}.
\end{equation}
If Equation~(\ref{EqnEqualPartials}) holds, the complex derivative
$\frac{\partial f}{\partial z}$ is defined to be the common value
of $\frac{\partial f}{\partial x}$ and $\frac{1}{i} \frac{\partial
f}{\partial y}$.  If Equation~(\ref{EqnEqualPartials}) does not
hold, then the function $f$ has no complex derivative.  For
example, complex conjugation $f(z) = \bar z$ fails
test~(\ref{EqnEqualPartials}) and the symbol $\frac{\partial\bar
z}{\partial z}$ has no meaning in this context.

An alternative definition of the derivative operator
\begin{equation}
  \frac{\partial}{\partial z}
   := \frac{1}{2}\left(\frac{\partial}{\partial x}
   + \frac{1}{i}\frac{\partial}{\partial y}\right)
\end{equation}
agrees with the traditional $\frac{\partial}{\partial z}$ when
applied to complex-differentiable functions, yet offers the
advantage of applying to the broader class of complex-valued
real-differentiable functions.  Any linear combination of
$\frac{\partial}{\partial x}$ and $\frac{1}{i}
\frac{\partial}{\partial y}$ would work, but the coefficients
$(\frac{1}{2}, \frac{1}{2})$ yield the desirable result that
\begin{equation}
  \frac{\partial \bar z}{\partial z}
  = \frac{1}{2}\left(\frac{\partial}{\partial x} + \frac{1}{i}
     \frac{\partial}{\partial y}\right) (x - i y)
  = 0.
\end{equation}
To take derivatives with respect to $\bar z$, define
\begin{equation}
  \frac{\partial}{\partial \bar z}
   := \frac{1}{2}\left(\frac{\partial}{\partial x}
   - \frac{1}{i}\frac{\partial}{\partial y}\right)
\end{equation}
and note that $\frac{\partial \bar z}{\partial \bar z} = 1$ and
$\frac{\partial z}{\partial \bar z} = 0$.  For both
$\frac{\partial}{\partial z}$ and $\frac{\partial}{\partial \bar
z}$ all the usual differentiation
rules (product rule, quotient rule, power rule) remain valid. \\

\section{The Laplace Operator}
\label{SectionLaplace}

Parameterize $\mathbb C^2 \approx R^4$ by two complex variables
$\alpha = x + i y$ and $\beta = z + i w$.  The form of the mixed
partials
\begin{equation}
    \frac{\partial^2}{\partial\alpha\;\partial\bar\alpha}
  = \frac{\partial^2}{\partial\bar\alpha\;\partial\alpha}
  = \frac{1}{4}\left(\frac{\partial^2}{\partial x^2} + \frac{\partial^2}{\partial y^2}\right)
\end{equation}
leads us immediately to the complex expression for the Laplace
operator
\begin{eqnarray}
\label{EqnLaplacian}
    \nabla^2
  &=& \frac{\partial^2}{\partial x^2} + \frac{\partial^2}{\partial y^2}
   +  \frac{\partial^2}{\partial z^2} + \frac{\partial^2}{\partial w^2}\nonumber\\
  &=& 4 \left(
      \frac{\partial^2}{\partial\alpha\;\partial\bar\alpha}
   +  \frac{\partial^2}{\partial\beta \;\partial\bar\beta }
   \right).
\end{eqnarray}

\section{Sibling Modes}
\label{SectionSiblingModes}

Each polynomial of homogeneous degree $k$ in $\alpha$ and $\beta$
alone (no $\bar\alpha$ or $\bar\beta$) generates a collection of
$k + 1$ sibling modes.  For example, the polynomial $\alpha^3 -
\beta^3$ generates the four sibling modes
\begin{eqnarray}
\label{EqnSampleSiblings}
  &\alpha^3 - \beta^3&\nonumber\\
  &\alpha^2\bar\beta + \beta^2\bar\alpha&\nonumber\\
  &\alpha\bar\beta^2 - \beta\bar\alpha^2&\nonumber\\
  &\bar\beta^3 + \bar\alpha^3&
\end{eqnarray}

Siblings are distinguished by differing {\it twist}.  Look what
happens as we trace the value of the original polynomial $\alpha^3
- \beta^3$ (the ``first sibling'') along an arbitrary but fixed
fiber $(\alpha, \beta) = (e^{i t} \alpha_0, e^{i t} \beta_0)$. The
polynomial evaluates to $e^{3 i t}(\alpha_0^3 - \beta_0^3)$ on the
fiber, so as $t$ runs from $0$ to $2\pi$ the polynomial's modulus
remains constant while its phase runs from $0$ to $6\pi$.  In
other words, the first sibling twists three times as we run along
a fiber. By contrast, the second sibling $\alpha^2\bar\beta +
\beta^2\bar\alpha = e^{(2 - 1) i t}(\alpha_0^2\bar\beta_0 +
\beta_0^2\bar\alpha_0)$ twists only once along a fiber;  the
regular variables $\alpha^2$ and $\beta^2$ contribute two positive
twists while the conjugated variables $\bar\beta$ and $\bar\alpha$
contribute a negative twist. Similarly the third sibling
$\alpha\bar\beta^2 - \beta\bar\alpha^2 = e^{(1 - 2) i
t}(\alpha_0\bar\beta_0^2 - \beta_0\bar\alpha_0^2)$ twists minus
once along a fiber and the last sibling $\bar\alpha^3 +
\bar\beta^3 = e^{-3 i t}(\bar\alpha_0^3 + \bar\beta_0^3)$ twists
minus three times.

Kramer~\cite{Kramer} generates siblings using what are commonly
called ``raising and lowering operators''.  Here we will define
the same thing but call them the {\it negative} and {\it positive
twist operators}
\begin{eqnarray}
\label{EqnTwistOperators}
  twist_{-} = -\bar\beta  \frac{\partial}{\partial\alpha}
            +  \bar\alpha \frac{\partial}{\partial\beta}\nonumber\\%
  twist_{+} =      -\beta \frac{\partial}{\partial\bar\alpha}
            +      \alpha \frac{\partial}{\partial\bar\beta }
\end{eqnarray}
to emphasize their geometrical significance, which will come in
handy in Section~\ref{SectionTwistAndInvariance}.  Algebraically
the twist operators' effect is clear:  $twist_{-}$ removes a
regular variable (via the partial derivatives) and replaces it
with a conjugated variable (via the multiplications) for a net
decrease of two twists, while $twist_{+}$ does the opposite.
Geometrically the factorization of a twist operator as a dot
product $(-\bar\beta,
\bar\alpha)\cdot(\frac{\partial}{\partial\alpha},
\frac{\partial}{\partial\beta })$ allows an interpretation as a
directional derivative orthogonal to the fiber. In any case, the
reader may verify that the positive and negative twist operators
take siblings up and down the
list~(\ref{EqnSampleSiblings}), modulo normalization. \\

\noindent{\bf Proposition \ref{SectionSiblingModes}.1.} {\it Modes
with different twists are orthogonal.}

\vskip0.25cm

\noindent{\it Proof.}  Let $f$ and $g$ be functions with twist $m$
and $n$, respectively.  When evaluating the integral $\langle
f,g\rangle = \int_{S^3}f\cdot\bar g$, note that along each fiber
$(\alpha, \beta) = (e^{i t} \alpha_0, e^{i t} \beta_0)$ the
integral restricts to
\begin{eqnarray}
  &\phantom{=}&
    \int_{t=0}^{2\pi}e^{i m t}f(\alpha_0,\beta_0)
                     e^{-i n t}\overline{g(\alpha_0,\beta_0)}\nonumber\\
  &=& f(\alpha_0,\beta_0)\overline{g(\alpha_0,\beta_0)}
      \int_{t=0}^{2\pi}e^{i (m - n) t}
\end{eqnarray}
which equals zero whenever $m \neq n$. $\blacksquare$ \\

\noindent{\bf Proposition \ref{SectionSiblingModes}.2.} {\it If a
symmetry group $\Gamma \subset SU(2)$ preserves a polynomial~$f$,
then it preserves all the siblings of $f$ as well.}

\vskip0.25cm

\noindent{\it Proof.}  Each symmetry $\gamma = ( \phantom{|}^z_w
\phantom{|}^{-\bar w}_{\phantom{-}\bar z} \phantom{|} ) \in
\Gamma$ defines a rotation $R_\gamma$ of function space via the
rule $[R_\gamma(f)](\alpha,\beta) = f[\gamma^{-1}(\alpha,\beta)]$
for all functions $f$ on $S^3 \subset \mathbb C^2$ and all points
$(\alpha,\beta) \in S^3$.  By assumption $R_\gamma(f) = f$.  A
quick calculation shows that $R_\gamma$ commutes with the positive
and negative twist operators $twist_\pm$ (alternatively,
interpreting the twist operators as directional derivatives
orthogonal to the fiber provides deeper geometrical insights, but
for sake of brevity we will not pursue that interpretation).  Thus
\begin{equation}
  R_\gamma[twist_\pm(f)] = twist_\pm[R_\gamma(f)] = twist_\pm(f).
\end{equation}
In other words, $\gamma$ preserves $twist_\pm(f)$ and therefore
preserves all siblings of $f$. $\blacksquare$

\section{Twist and $\Gamma$-Invariance}
\label{SectionTwistAndInvariance}

The twist concept sheds some light on why symmetric polynomials
turn out to be invariant under symmetries of their roots in some
cases (Example~\ref{SectionSymmetricPolynomials}.1) but not others
(Example~\ref{SectionSymmetricPolynomials}.2). Geometrical
considerations often reveal at a glance that a given mode cannot
possibly be invariant, eliminating the need to test the invariance
explicitly as we did in Section~\ref{SectionSymmetricPolynomials}.
\\

\noindent{\bf Example \ref{SectionTwistAndInvariance}.1.} {\it The
octahedral group $O$.}  An order-2 rotation about the midpoint of
any of the octahedron's edges lifts to an order-4 corkscrew motion
of $S^3$ (Section~\ref{SectionGroupConstruction}). That order-4
corkscrew motion fixes (setwise) the fiber in $S^3$ corresponding
to the fixed edge midpoint on $S^2$.  Therefore any
$O^*$-invariant function must have period~4 along the fiber,
meaning its twist must be a multiple of~4.  Thus the degree-6
polynomial constructed from the octahedron's six vertices, which
has twist~6 along every fiber, cannot possibly be $O^*$-invariant.
\\

\noindent{\bf Example \ref{SectionTwistAndInvariance}.2.} {\it The
tetrahedral group $T$.}  As in the preceding example, an order-2
rotation about the midpoint of any edge lifts to an order-4
corkscrew motion of $S^3$, implying that any $T^*$-invariant
function must have period~4 along the fibers corresponding to
fixed points.  This would seem to exclude the degree-6 polynomial
constructed from the tetrahedron's six edge midpoints, but in fact
that degree-6 polynomial -- by construction! -- is identically
zero along the fixed fibers. Thus the degree-6 polynomial could
well be (and indeed is) invariant under the period-4 corkscrew
motion, and in fact turns out to be invariant under all of $T^*$.

On the other hand, the degree-4 polynomial constructed from the
tetrahedron's face centers cannot possibly be $T^*$-invariant,
because along a fiber it's not periodic of order~6, as the order-3
rotation about any vertex would require, and is nonzero along half
the fixed fibers.
\\

\noindent{\bf Example \ref{SectionTwistAndInvariance}.3.} {\it The
icosahedral group $I$.}  Section~\ref{SectionModesS3I} will
construct the degree~12 polynomial $I_{12}$ corresponding to the
icosahedron's 12~vertices, the degree~30 polynomial $I_{30}$
corresponding to the 30~edge midpoints, and the degree~20
polynomial $I_{20}$ corresponding to the 20~faces
(Equation~(\ref{EqnModesIcos})).  As it turns out, luck has smiled
on the icosahedral group.
\begin{itemize}
  \item The 5-fold rotation about each vertex requires that an
        invariant polynomial have period~10 along each fiber
        corresponding to a fixed point.  The polynomials $I_{20}$
        and $I_{30}$ satisfy this criterion automatically because
        they have twist 20 and 30, respectively.  While the twist
        of $I_{12}$ is not a multiple of 10, its zeros fall
        at the fixed points of the rotations, so it is identically
        zero along the fibers in question.
  \item The 3-fold rotation about each face center requires
        period-6 invariance along the fixed fibers.  The polynomials
        $I_{12}$ and $I_{30}$, whose twists are multiples of $6$,
        satisfy this automatically.  For $I_{20}$ the degree is
        wrong, but luckily it's identically zero along the fixed
        fibers, so it is safe too.
  \item The 2-fold rotation about each edge midpoint requires
        4-fold periodicity along the fiber, which $I_{12}$ and $I_{20}$
        satisfy by virtue of their twist, and $I_{30}$ satisfies
        because it's identically zero on the fixed fibers.
\end{itemize}
Thus we find no geometrical obstructions to the $I^*$-invariance
of $I_{12}$, $I_{20}$ and $I_{30}$ and the method of
Section~\ref{SectionSymmetricPolynomials} confirms all three to be
fully $I^*$-invariant.
\\

\section{Harmonicity of Modes}
\label{SectionHarmonicity}

By applying the complex form of the Laplace
operator~(\ref{EqnLaplacian}) to the sample
siblings~(\ref{EqnSampleSiblings}), the reader may quickly verify
that all are harmonic.  The following propositions show that this
will always be
the case.\\

\noindent{\bf Proposition \ref{SectionHarmonicity}.1.} {\it Every
polynomial in $\alpha$ and $\beta$ alone (no $\bar\alpha$ or
$\bar\beta$) is harmonic.}

\vskip0.25cm

\noindent{\it Proof.}  If $f$ is a polynomial in $\alpha$ and
$\beta$ alone, then $\frac{\partial f }{\partial\bar\alpha} =
\frac{\partial f }{\partial\bar\beta} = 0$ and thus $\nabla^2 f =
4 \left( \frac{\partial^2 f }{\partial\alpha\;\partial\bar\alpha}
+ \frac{\partial^2 f}{\partial\beta \;\partial\bar\beta } \right)
= 0$.
$\blacksquare$ \\

\noindent{\bf Proposition \ref{SectionHarmonicity}.2.} {\it If a
polynomial is harmonic then so are all its siblings.}

\vskip0.25cm

\noindent{\it Proof.}  The operators $\nabla^2$ and $twist_{-}$
commute, so if $\nabla^2 f = 0$, then $\nabla^2 (twist_{-}\,f) =
twist_{-} (\nabla^2 f) = 0$ as well, and similarly for
$twist_{+}$.
$\blacksquare$ \\

\noindent Every harmonic function on $\mathbb C^2 \approx \mathbb
R^4$ restricts to an eigenmode of the Laplace operator
$\nabla^2_{S^3}$ on $S^3$, because
\begin{equation}
    \nabla^2_{S^3}
  = {\bf r}^2\nabla^2_{\mathbb R^4}
  - ({\bf r}\cdot\nabla)({\bf r}\cdot\nabla)
  - 2({\bf r}\cdot\nabla)
\end{equation}
which for a polynomial $Q$ of homogeneous degree $k$ simplifies to
\begin{equation}
    \nabla^2_{S^3}Q
  = {\bf r}^2\nabla^2_{\mathbb R^4}Q
  - k(k+2)Q
\end{equation}
and when $Q$ is harmonic ($\nabla^2_{\mathbb R^4}Q = 0$)
simplifies further to
\begin{equation}
    \nabla^2_{S^3}Q
  =  - k(k+2)Q.
\end{equation}

\section{Modes of the Poincar\'e Dodecahedral Space $S^3/I^*$}
\label{SectionModesS3I}

For each wavenumber $k$ we will define zero or more {\it base
modes}.  These base modes, together with their siblings
(Section~\ref{SectionSiblingModes}), will form a complete
eigenbasis for the Poincar\'e dodecahedral space $S^3/I^*$ for the
given $k$.

For $k = 0$ the base mode is the constant polynomial $1$, which of
course has no siblings.

For $k = 12$, $30$ or $20$, construct symmetric polynomials as in
Section~\ref{SectionSymmetricPolynomials}, taking $P$ to be the
set of vertices, the set of edge midpoints or the set of face
centers, respectively, of a regular icosahedron, yielding the base
modes
\begin{eqnarray}
\label{EqnModesIcos}
  I_{12} &=& \alpha^{11} \beta + 11 \alpha^6 \beta^6 - \alpha\beta^{11}\nonumber\\
  I_{30} &=& \alpha^{30} + 522 \alpha^{25} \beta^5 - 10005 \alpha^{20} \beta^{10}%
    - 10005 \alpha^{10} \beta^{20} - 522 \alpha^5 \beta^{25} + \beta^{30}\nonumber\\
  I_{20} &=& \alpha^{20} - 228 \alpha^{15} \beta^5 + 494 \alpha^{10} \beta^{10}%
    + 228 \alpha^5 \beta^{15} + \beta^{20}.
\end{eqnarray}
Still using the method of
Section~\ref{SectionSymmetricPolynomials}, one easily verifies
that these three polynomials are fully invariant under the binary
icosahedral group $I^* \subset SU(2)$.  These three base modes,
together with their siblings, form complete eigenbases of
dimension $13$, $31$ and $21$, respectively.  Historical note:
Klein's monograph~\cite{Klein} presents these polynomials and
notes their invariance under the ``group of substitutions'', but
never considers the quotient space $S^3/I^*$. Clifford-Klein space
forms weren't conceived until a few years later.

For wavenumbers $k \neq 0, 12, 20, 30$, form the base modes as
products of $I_{12}$, $I_{20}$ and $I_{30}$.  For example, for $k
= 32$ there is exactly one base mode, namely $I_{12}I_{20}$, which
together with its siblings forms a basis of dimension $33$.  For
$k = 34$ there are no base modes, because $34$ cannot be expressed
as a sum of 12's, 20's and 30's.  All products of $I_{12}$,
$I_{20}$ and $I_{30}$ are harmonic by Proposition
\ref{SectionHarmonicity}.1 and clearly all are invariant under the
binary icosahedral group $I^*$.

The wavenumber $k = 60$ is especially interesting because $60$ may
be written as a sum of 12's, 20's and 30's in more than one way,
namely $60 = 5\cdot 12 = 3\cdot 20 = 2\cdot 30$.  This gives three
candidates for the base mode, namely $I_{12}^5$, $I_{20}^3$ and
$I_{30}^2$.  However, a simple computation reveals a linear
dependence among the three modes, so any two will span the same
2-dimensional function space.  For sake of concreteness, let us
choose $I_{12}^5$ and $I_{20}^3$, ignoring $I_{30}^2$.  (In
cosmological applications one typically orthogonalizes this basis,
but for present purposes this isn't necessary.) The two base modes
$I_{12}^5$ and $I_{20}^3$, together with their siblings, form a
complete eigenbasis of dimension $61 + 61 = 122$. \\

\noindent{\bf Definition \ref{SectionModesS3I}.1.}  The {\it
modulus function} $\mathrm{mod}(m,n)$ gives the remainder of $m$
upon integer division by $n$.  For example, $\mathrm{mod}(17,5) =
2$. \\

\noindent{\bf Definition \ref{SectionModesS3I}.2.}  The {\it floor
function} $\lfloor x\rfloor$ gives the integer part of its
argument $x$. For example, $\lfloor\frac{17}{5}\rfloor = 3$. \\

The Poincar\'e dodecahedral space has no modes for odd $k$,
because odd modes fail to be invariant under the antipodal map
$-id \in I^*$.  To enumerate the base modes for arbitrary even
$k$, it is convenient to introduce the notation $\hat k = k/2$.
Rather than asking how $k$ may be written as a sum $a\cdot 12 +
b\cdot 20 + c\cdot 30$, we'll ask how $\hat k$ may be written as a
sum
\begin{eqnarray}
\label{EqnKHat}
  \hat k &=& a\cdot 6 + b\cdot 10 + c\cdot 15 \nonumber\\
         &=& a\cdot (2 \cdot 3) + b\cdot (2 \cdot 5) + c\cdot (3 \cdot 5).
\end{eqnarray}
First consider $\hat k$ modulo~5.  From (\ref{EqnKHat}) it is
clear that $a \equiv \hat k \;(\mathrm{mod}\;5)$, so in particular
$a \geq \mathrm{mod}(\hat k,5) = \hat k - 5\lfloor \frac{\hat
k}{5}\rfloor$, with notation as in
Definitions~\ref{SectionModesS3I}.1 and \ref{SectionModesS3I}.2.
Similarly, $b \geq \mathrm{mod}(\hat k,3) = \hat k - 3\lfloor
\frac{\hat k}{3}\rfloor$ and $c \geq \mathrm{mod}(\hat k,2) = \hat
k - 2\lfloor \frac{\hat k}{2}\rfloor$. Subtracting off the minimum
number of 6's, 10's and 15's leaves
\begin{eqnarray}
\label{EqnKHatRemainder}
  &\phantom{=}& \hat k
    - \left(\hat k - 5\left\lfloor\frac{\hat k}{5}\right\rfloor\right)  6
    - \left(\hat k - 3\left\lfloor\frac{\hat k}{3}\right\rfloor\right) 10
    - \left(\hat k - 2\left\lfloor\frac{\hat k}{2}\right\rfloor\right) 15\nonumber\\
  &=& 30
    \left(
        \left\lfloor\frac{\hat k}{5}\right\rfloor
      + \left\lfloor\frac{\hat k}{3}\right\rfloor
      + \left\lfloor\frac{\hat k}{2}\right\rfloor
      - \hat k
    \right)
\end{eqnarray}
The quantity
\begin{equation}
  \tilde k =
        \left\lfloor\frac{\hat k}{5}\right\rfloor
      + \left\lfloor\frac{\hat k}{3}\right\rfloor
      + \left\lfloor\frac{\hat k}{2}\right\rfloor
      - \hat k
\end{equation}
tells how many groups of 30 remain after we have subtracted off
the minimum 6's, 10's and 15's.

If $\tilde k = -1$, then no decomposition $\hat k = a\cdot 6 +
b\cdot 10 + c\cdot 15$ is possible (with nonnegative integer
coefficients) and the Poincar\'e dodecahedral space has no modes
of the given order $k$. If $\tilde k = 0$, then the decomposition
is unique, consisting of exactly the required minimum number of
6's, 10's and 15's, and the mode $I_{12}^{a} I_{20}^{b}
I_{30}^{c}$ along with its siblings spans an eigenspace of
dimension $k + 1$.  If $\tilde k > 0$ then there will be $\tilde
k$ factors of degree $60$, which may be apportioned between
$I_{12}^5$ and $I_{20}^3$ in $\tilde k + 1$ different ways
(ignoring $I_{30}^2$ because it is a linear combination of
$I_{12}^5$ and $I_{20}^3$). \\

\noindent{\bf Example \ref{SectionModesS3I}.3.}  When $k = 180$,
then $\hat k = 90$ and $\tilde k = 3$.  The $\tilde k+1=4$ base
modes
\begin{equation}
    (I_{12}^5)^3,\quad
    (I_{12}^5)^2(I_{20}^3),\quad
    (I_{12}^5)(I_{20}^3)^2,\quad
    (I_{20}^3)^3,
\end{equation}
each with its $k+1=181$ siblings, span an eigenspace of dimension
$4\cdot 181 = 724$. \\

\noindent{\bf Example \ref{SectionModesS3I}.4.}  When $k = 234$,
then $\hat k = 117 = 3\cdot 30 + 2\cdot 6 + 1\cdot 15$ and there
are still the same $\tilde k = 3$ factors of degree $60$ as in
Example~\ref{SectionModesS3I}.1, but now they are supplemented by
required minimum factors of degree $12$ and $30$. Specifically,
$I_{12}$ must appear with exponent $a = \mathrm{mod}(\hat k,5) =
2$ and $I_{30}$ must appear with exponent $c = \mathrm{mod}(\hat
k,2) = 1$.  No extra factor of $I_{20}$ is required because $b =
\mathrm{mod}(\hat k,3) = 0$.  Thus for $k = 234$ the base modes
are
\begin{equation}
  (I_{12}^2I_{30})(I_{12}^5)^3,\;
  (I_{12}^2I_{30})(I_{12}^5)^2(I_{20}^3),\;
  (I_{12}^2I_{30})(I_{12}^5)(I_{20}^3)^2,\;
  (I_{12}^2I_{30})(I_{20}^3)^3.
\end{equation}
Together with their siblings they span an eigenspace of dimension
$4\cdot 235 = 940$. \\

\noindent{\bf Theorem \ref{SectionModesS3I}.5.} {\it For each even
$k$, if we let $i$ range from $0$ to
\begin{equation}
\label{EqnDimensionI}
  \tilde k =
      \left\lfloor\frac{k/2}{2}\right\rfloor
    + \left\lfloor\frac{k/2}{3}\right\rfloor
    + \left\lfloor\frac{k/2}{5}\right\rfloor
    - \frac{k}{2},
\end{equation}
then the $\tilde k+1$ base polynomials
\begin{equation}
\begin{array}{c}
    (\alpha^{11} \beta + 11 \alpha^6 \beta^6 - \alpha \beta^{11})^{5i + \mathrm{mod}(\frac{k}{2},5)} \\
    \times\;(\alpha^{20} - 228 \alpha^{15} \beta^5 + 494 \alpha^{10} \beta^{10}
          + 228 \alpha^5 \beta^{15} + \beta^{20})^{3(\tilde k - i) + \mathrm{mod}(\frac{k}{2},3)} \\
    \times\;(\alpha^{30} + 522 \alpha^{25} \beta^5 - 10005 \alpha^{20} \beta^{10}
        - 10005 \alpha^{10} \beta^{20} - 522 \alpha^5 \beta^{25} + \beta^{30})^{\mathrm{mod}(\frac{k}{2},2)},
\end{array}
\end{equation}
each with its $k+1$ siblings (defined in
Section~\ref{SectionSiblingModes}), form a basis of dimension
$(\tilde k + 1)(k + 1)$ for the eigenmodes of the Poincar\'e
dodecahedral space $S^3/I^*$. There are no modes for odd $k$.}

\vskip0.25cm

\noindent{\it Proof.}  To see that the $\tilde k + 1$ base modes
are linearly independent, note that they all have different
leading terms in the variable $\alpha$, because $\alpha$ occurs
with maximal degree $55$ in $I_{12}^5$ but maximal degree $60$ in
$I_{20}^3$ (Equation~(\ref{EqnModesIcos})).  Allowing for the base
modes' siblings, the same reasoning applies, but instead of
measuring powers of $\alpha$ alone, measure combined powers of
$\alpha$ and $\bar\beta$, a quantity which the twist
operator~(\ref{EqnTwistOperators}) preserves.  This proves that
linear dependencies cannot occur among the siblings of different
base modes. Similar reasoning applies within each family of
siblings, where differing twist expresses itself algebraically as
differing powers of the regular variables $\alpha$ and $\beta$
relative to the conjugated variables $\bar\alpha$ and $\bar\beta$.
Therefore the full collection of base modes and siblings is
linearly independent.  Technical note:  This line of reasoning
relies on the unique expression of a polynomial in terms of
$\alpha$, $\bar\alpha$, $\beta$ and $\bar\beta$, which follows
easily from its uniqueness in terms of the underlying real
variables $x$, $y$, $z$ and $w$.

The $\tilde k+1$ base modes, each with $k + 1$ siblings (including
itself), span a space of dimension
\begin{equation}
  (\tilde k + 1)(k + 1) =
  \left(\left\lfloor\frac{k/2}{5}\right\rfloor
      + \left\lfloor\frac{k/2}{3}\right\rfloor
      + \left\lfloor\frac{k/2}{2}\right\rfloor
      - \frac{k}{2} + 1\right)
  (k + 1).
\end{equation}
This agrees with Ikeda's formula for the dimension of the
eigenspace~\cite[Thm.~4.6]{Ikeda}, proving that our modes span the
full eigenspace of $S^3/I^*$.  Note:  Ikeda computes the spectrum
using a dimension counting argument, without constructing the
modes, so no further comparison is possible.

As we saw earlier, there are no modes for odd $k$ because such
modes are not invariant under the antipodal map $-id \in I^*$.
$\blacksquare$ \\

Reasoning similar to that in the preceding proof shows that for a
given base mode, sibling $i$ must be the complex conjugate of
sibling $k - i$, for $i=0\dots k$, up to normalization.  In
practice the normalization is chosen so that corresponding
siblings are exactly conjugate to each other, to facilitate easy
extraction of the real-valued modes.  The middle sibling ($i =
\frac{k}{2}$) is of course its own conjugate up to a factor of
$\pm 1$.

\section{Modes of the Binary Octahedral Space $S^3/O^*$}
\label{SectionModesS3O}

Let us apply the method of previous sections to write down the
eigenmodes of the binary octahedral space $S^3/O^*$.  Consider in
turn the octahedron's vertices, its edge midpoints and its faces
centers, which according to the method of
Section~\ref{SectionSymmetricPolynomials} yield the tentative base
modes
\begin{eqnarray}
\label{EqnModesOctTentative}
  O_{6\phantom{0}}^\prime &=& \alpha^5\beta - \alpha\beta^5 \nonumber\\
  O_{12}^\prime           &=& \alpha^{12} - 33\alpha^8\beta^4 - 33\alpha^4\beta^8 + \beta^{12} \nonumber\\
  O_{8\phantom{0}}^\prime &=& \alpha^8 + 14\alpha^4\beta^4 + \beta^8.
\end{eqnarray}
The catch is that these base modes are not all $O^*$-invariant, as
Example~\ref{SectionTwistAndInvariance}.1 illustrated.  More
precisely, a quarter turn rotation about any vertex of the
octahedron lifts to an order~8 Clifford translation of $S^3$ that
preserves $O_{8}^\prime$ while sending $O_{6}^\prime$ and
$O_{12}^\prime$ to their negatives. Thus neither $O_{6}^\prime$
nor $O_{12}^\prime$ defines a mode of $S^3/O^*$ on its own, but
their product $O_{6}^\prime O_{12}^\prime$ and their squares
$O_{6}^{\prime 2}$ and $O_{12}^{\prime 2}$ do.  We may safely
ignore $O_{12}^{\prime 2}$ because the mode it represents can be
recovered later as a linear combination $O_{12}^{\prime 2} =
O_{8}^{\prime 3} - 108\,O_{6}^{\prime 4}$.  Thus we take as our
base modes
\begin{equation}
  \begin{array}{lclcl}
    O_{ 8} &=& O_{8}^\prime &=& \alpha^8 + 14\alpha^4\beta^4 + \beta^8 \\
    O_{12} &=& O_{6}^{\prime\;2} &=& \alpha^{10}\beta^2 - 2\alpha^6\beta^6 + \alpha^2\beta^{10} \\
    O_{18} &=& O_{6}^\prime O_{12}^\prime &=& \alpha^{17}\beta - 34\alpha^{13}\beta^5 + 34\alpha^5\beta^{13} - \alpha\beta^{17}.\\
  \end{array}
\end{equation}

Now proceed exactly as in Section~\ref{SectionModesS3I}, but using
$\{O_{8}, O_{12}, O_{18}\}$ instead of $\{I_{12}, I_{20},
I_{30}\}$.  That is, for each even wavenumber $k$, ask how many
base modes (if any) may be constructed as products of the
$\{O_{8}, O_{12}, O_{18}\}$.  Still imitating
Section~\ref{SectionModesS3I}, define $\hat k = k/2$ and ask how
$\hat k$ may be written as a sum
\begin{equation}
\label{EqnOctKHat}
  \hat k = a\cdot 4 + b\cdot 6 + c\cdot 9.
\end{equation}
One sees immediately that $a \equiv \hat k \;(\mathrm{mod}\;3)$,
$c \equiv \hat k \;(\mathrm{mod}\;2)$, and $2b \equiv \hat k - c
\;(\mathrm{mod}\;4)$.  The last equivalence reduces to $b \equiv
\frac{\hat k - c}{2} \;(\mathrm{mod}\;2)$.  Furthermore one may
insist that $c$ be $0$ or $1$, because any occurrence of
$O_{18}^2$ may be replaced by $O_{18}^2 = O_{6}^{\prime 2}
O_{12}^{\prime 2} = O_{6}^{\prime 2} (O_{8}^{\prime 3} -
108\,O_{6}^{\prime 4}) = O_{12} (O_{8}^3 - 108\,O_{12}^2)$.
Setting $c = \mathrm{mod}(\hat k,2)$ transforms the restriction on
$b$ to $b \equiv \frac{\mathrm{mod}(\hat k,4) - \mathrm{mod}(\hat
k,2)}{2} \;(\mathrm{mod}\;2)$.  Subtracting from $\hat k$ the
minimum multiples of 4, 6 and 9 leaves
\begin{eqnarray}
\label{EqnOctKHatRemainder}
  &\hat k
    - \left(\hat k - 3\left\lfloor\frac{\hat k}{3}\right\rfloor\right) 4
    - \left(\frac
      {
        \left(\hat k - 4\left\lfloor\frac{\hat k}{4}\right\rfloor\right)
        -
        \left(\hat k - 2\left\lfloor\frac{\hat k}{2}\right\rfloor\right)
      }
      {2}\right)
      6
    - \left(\hat k - 2\left\lfloor\frac{\hat k}{2}\right\rfloor\right) 9&\nonumber\\
  &= 12
    \left(
        \left\lfloor\frac{\hat k}{2}\right\rfloor
      + \left\lfloor\frac{\hat k}{3}\right\rfloor
      + \left\lfloor\frac{\hat k}{4}\right\rfloor
      - \hat k
    \right)&
\end{eqnarray}
and so the quantity
\begin{equation}
  \tilde k =
        \left\lfloor\frac{\hat k}{2}\right\rfloor
      + \left\lfloor\frac{\hat k}{3}\right\rfloor
      + \left\lfloor\frac{\hat k}{4}\right\rfloor
      - \hat k
\end{equation}
tells how many groups of 12 remain.  Each group of 12 may be
realized as either $O_{8}^3$ or $O_{12}^2$. \\

\noindent{\bf Theorem \ref{SectionModesS3O}.1.} {\it For each even
$k$, if we let $i$ range from $0$ to
\begin{equation}
\label{EqnDimensionO}
  \tilde k =
      \left\lfloor\frac{k/2}{2}\right\rfloor
    + \left\lfloor\frac{k/2}{3}\right\rfloor
    + \left\lfloor\frac{k/2}{4}\right\rfloor
    - \frac{k}{2},
\end{equation}
then the $\tilde k+1$ base polynomials
\begin{equation}
\begin{array}{c}
    (\alpha^8 + 14\alpha^4\beta^4 + \beta^8)^{3i + \mathrm{mod}(\frac{k}{2},3)} \\
    \times\;(\alpha^{10}\beta^2 - 2\alpha^6\beta^6 + \alpha^2\beta^{10})
         ^{2(\tilde k - i) + \left(\mathrm{mod}(\frac{k}{2},4) - \mathrm{mod}(\frac{k}{2},2)\right)/2} \\
    \times\;(\alpha^{17}\beta - 34\alpha^{13}\beta^5 + 34\alpha^5\beta^{13} -
         \alpha\beta^{17})^{\mathrm{mod}(\frac{k}{2},2)},\\
\end{array}
\end{equation}
each with its $k+1$ siblings (defined in
Section~\ref{SectionSiblingModes}), form a basis of dimension
$(\tilde k + 1)(k + 1)$ for the eigenmodes of $S^3/O^*$.  There
are no modes for odd $k$.}

\vskip0.25cm

\noindent{\it Proof.}  Same as the proof of
Theorem~\ref{SectionModesS3I}.5.  $\blacksquare$ \\

\section{Modes of the Binary Tetrahedral Space $S^3/T^*$}
\label{SectionModesS3T}

Let us apply the methods of Sections~\ref{SectionModesS3I} and
\ref{SectionModesS3O} to write down the eigenmodes of the binary
tetrahedral space $S^3/T^*$.  The tetrahedron's vertices, edge
midpoints and face centers provide tentative base modes
\begin{eqnarray}
  T_{4a}^\prime           &=& \alpha^4 + 2 i \alpha^2\beta^2 + \beta^4 \nonumber\\
  T_{6\phantom{a}}^\prime &=& \alpha^5\beta - \alpha\beta^5 \nonumber\\
  T_{4b}^\prime           &=& \alpha^4 - 2 i \alpha^2\beta^2 + \beta^4.
\end{eqnarray}
A one-third rotation about any vertex of the tetrahedron lifts to
an order~6 Clifford translation of $S^3$ that preserves
$T_{6}^\prime$ but sends $T_{4a}^\prime \mapsto
\bar\xi\;T_{4a}^\prime$ and $T_{4b}^\prime \mapsto
\xi\;T_{4b}^\prime$, where $\xi = \frac{1}{2} +
i\frac{\sqrt{3}}{2}$ is a third root of unity. Therefore the group
$T^*$ preserves the modes
\begin{equation}
\label{EqnModesTet}
  \begin{array}{lclcl}
    T_{ 6} &=& T_{6}^\prime
      &=& \alpha^5\beta - \alpha\beta^5
      \\
    T_{ 8} &=& T_{4a}^\prime T_{4b}^\prime
      &=& \alpha^8 + 14\alpha^4\beta^4 + \beta^8
      \\
    T_{12} &=& (T_{4a}^{\prime\;3} + T_{4b}^{\prime\;3})/2
      &=& \alpha^{12} - 33\alpha^8\beta^4 - 33\alpha^4\beta^8 + \beta^{12}.
      \\
  \end{array}
\end{equation}
We may ignore the complementary order~12 mode $(T_{4a}^{\prime\;3}
- T_{4b}^{\prime\;3})/2$ because it is a multiple of $T_{6}^2$.
Note that the tetrahedral modes~(\ref{EqnModesTet}) coincide
exactly with the tentative octahedral
modes~(\ref{EqnModesOctTentative}), reflecting the close ties
between $O^*$ and its index~2 subgroup $T^*$.

Let us construct the modes of $S^3/T^*$ as in the preceding
sections, by taking products of $\{T_{6},T_{8},T_{12}\}$ to
establish the base modes and then taking all their siblings.  With
$\hat k = k/2$ as before, ask how $\hat k$ may be written as a sum
\begin{equation}
\label{EqnTetKHat}
  \hat k = a\cdot 3 + b\cdot 4 + c\cdot 6.
\end{equation}
One sees immediately that $a \equiv \hat k \;(\mathrm{mod}\;2)$
and $b \equiv \hat k \;(\mathrm{mod}\;3)$.  Subtracting from $\hat
k$ the minimum multiples of 3 and 4 leaves
\begin{eqnarray}
\label{EqnTetKHatRemainder}
  &\hat k
    - \left(\hat k - 2\left\lfloor\frac{\hat k}{2}\right\rfloor\right) 3
    - \left(\hat k - 3\left\lfloor\frac{\hat k}{3}\right\rfloor\right) 4&\nonumber\\
  &= 6
    \left(
        \left\lfloor\frac{\hat k}{2}\right\rfloor
      + \left\lfloor\frac{\hat k}{3}\right\rfloor
      + \left\lfloor\frac{\hat k}{3}\right\rfloor
      - \hat k
    \right)&
\end{eqnarray}
and so the quantity
\begin{equation}
  \tilde k =
        \left\lfloor\frac{\hat k}{2}\right\rfloor
      + \left\lfloor\frac{\hat k}{3}\right\rfloor
      + \left\lfloor\frac{\hat k}{3}\right\rfloor
      - \hat k
\end{equation}
tells how many groups of 6 remain.  Each group of 6 may be
realized as either $T_{6}^2$ or $T_{12}$.  We may safely restrict
$T_{8}$ to the minimum $b = 0$, $1$ or $2$ occurrences (according
to the value of $k$~mod~3), because additional occurrences may be
swapped out via the dependency $T_{8}^3 = T_{12}^2 + 108 T_{6}^4$.
\\

\noindent{\bf Theorem \ref{SectionModesS3T}.1.} {\it For each even
$k$, if we let $i$ range from $0$ to
\begin{equation}
\label{EqnDimensionT}
  \tilde k =
      \left\lfloor\frac{k/2}{2}\right\rfloor
    + \left\lfloor\frac{k/2}{3}\right\rfloor
    + \left\lfloor\frac{k/2}{3}\right\rfloor
    - \frac{k}{2},
\end{equation}
then the $\tilde k+1$ base polynomials
\begin{equation}
\begin{array}{c}
    (\alpha^5\beta - \alpha\beta^5)^{2i + \mathrm{mod}(\frac{k}{2},2)}\;
    (\alpha^8 + 14\alpha^4\beta^4 + \beta^8)^{\mathrm{mod}(\frac{k}{2},3)} \\
    \times\;(\alpha^{12} - 33\alpha^8\beta^4 - 33\alpha^4\beta^8 + \beta^{12})^{\tilde k - i},\\
\end{array}
\end{equation}
each with its $k+1$ siblings (defined in
Section~\ref{SectionSiblingModes}), form a basis of dimension
$(\tilde k + 1)(k + 1)$ for the eigenmodes of $S^3/T^*$.  There
are no modes for odd $k$.}

\vskip0.25cm

\noindent{\it Proof.}  Same as the proof of
Theorems~\ref{SectionModesS3I}.5 and \ref{SectionModesS3O}.1.
$\blacksquare$ \\

\section{Modes of the Binary Dihedral Spaces $S^3/D_m^*$}
\label{SectionModesS3D}

Still following the method of the preceding three sections, let us
write down the eigenmodes of the binary dihedral space
$S^3/D_m^*$, where $D_m^*$ is the binary dihedral group of order
$4m$.  The dihedron's vertices (sitting at $m^{th}$ roots of~$1$),
edge midpoints (at $m^{th}$ roots of $-1$) and face centers (at
$\{0,\infty\}$) provide the tentative base modes
\begin{eqnarray}
\label{EqnDihedralBaseModesTentative}
  D_v^\prime &=& \alpha^m - \beta^m \nonumber\\
  D_e^\prime &=& \alpha^m + \beta^m \nonumber\\
  D_f^\prime &=& \alpha\beta.
\end{eqnarray}
A half turn about the vertex at $1$ lifts to the order~4 Clifford
translation $(\alpha,\beta) \mapsto ({}^0_i
{\,}^i_0)({}^\alpha_\beta) = (i\beta,i\alpha)$, which sends
$D_v^\prime \mapsto -i^m D_v^\prime$, $D_e^\prime \mapsto +i^m
D_e^\prime$, and $D_f^\prime \mapsto -D_f^\prime$.  Similarly a
half turn about the edge midpoint at $e^{\frac{i\pi}{m}}$ lifts to
the order~4 Clifford translation $(\alpha,\beta) \mapsto
\left({}^{\;\; 0}_{i e^{\frac{-i\pi}{m}}} {\,}^{i
e^{\frac{i\pi}{m}}}_{\;\; 0}\right)({}^\alpha_\beta) = (i
e^{\frac{i\pi}{m}} \beta, i e^{\frac{-i\pi}{m}} \alpha)$, which
sends $D_v^\prime \mapsto +i^m D_v^\prime$, $D_e^\prime \mapsto
-i^m D_e^\prime$, and $D_f^\prime \mapsto -D_f^\prime$.  These two
half turns generate the full group.

Keeping in mind that the product $D_v^\prime D_e^\prime =
\alpha^{2m} - \beta^{2m}$ is invariant when $m$ is odd, while
$(D_e^\prime)^2 - 2(D_f^\prime)^m = \alpha^{2m} + \beta^{2m}$ is
invariant when $m$ is even, we may combine the tentative base
modes~(\ref{EqnDihedralBaseModesTentative}) into the final base
modes
\begin{equation}
\label{EqnDihedralBaseModes}
\begin{array}{lcl}
  D_{4}    & = & \alpha^2 \beta^2 \\
  D_{2m}   & = & \alpha^{2m} \pm \beta^{2m} \\
  D_{2m+2} & = & \alpha\beta(\alpha^{2m} \mp \beta^{2m}) \\
\end{array}
\end{equation}
with the upper sign choices for $m$ even and the lower sign
choices for $m$ odd.  In both cases the modes are fully invariant
under the two half-turns and therefore under all of~$D_m^*$.

For each wavenumber $k$ ask how many independent base modes may be
constructed as a product
\begin{equation}
  D_{4}^{a} D_{2m}^{b} D_{2m+2}^{c}.
\end{equation}
We may safely insist that $c$ be $0$ or $1$, because the square
$D_{2m+2}^2$ may be replaced by $D_{2m+2}^2 = D_{4}(D_{2m}^2 \mp 4
D_{4}^m)$.  To facilitate an inductive approach, let $B_k$ denote
the set of modes constructed for wavenumber $k$.  For $k = 0$ the
only mode is the polynomial $1$, so $B_0 = \{ 1 \}$. For $k = 2$
there are no base modes when $m > 1$, so $B_2 = \{\}$ is empty,
but there is a base mode $B_2 = \{\alpha^{2} - \beta^{2}\}$ when
$m = 1$.  Now proceed inductively.  Assuming we've already
constructed the set $B_{k-4}$, let $B_{k}$ consist of the
polynomials in $B_{k-4}$ each multiplied by an additional factor
of $D_{4}$, along with an additional mode
$D_{2m}^{(\frac{k}{2m})}$ if $2m$ divides $k$ or
$D_{2m}^{(\frac{k-(2m+2)}{2m})} D_{2m+2}$ if $2m$ divides $k-2$
(or both in the special case $m = 1$).  To prove that the
polynomials in $B_{k}$ are linearly independent, note that each
one has a different highest power of $\alpha$.

A few concrete examples illustrate the iterative --- and
ultimately very simple -- nature of the construction.  The reader
is encouraged to carry out the construction, as defined in the
preceding paragraph, for the cases $m = 1,2,3$ and compare his or
her results to those given in the following three examples. \\

\noindent{\bf Example \ref{SectionModesS3D}.1.} In the seemingly
(but not really) exceptional case $m = 1$ the construction
proceeds as follows:
\begin{equation}
\begin{array}{lcl}
  B_{0}  & = & \{ 1 \} \\
  B_{2}  & = & \{ \alpha^{2} - \beta^{2} \} \\
  B_{4}  & = & \{ \alpha^2\beta^2,\;
                  (\alpha^2 - \beta^2)^2,\;
                  \alpha\beta(\alpha^2 + \beta^2) \} \\
  B_{6}  & = & \{ \alpha^2\beta^2(\alpha^{2} - \beta^{2}),\;
                  (\alpha^2 - \beta^2)^3,\;
                  (\alpha^2 - \beta^2)\alpha\beta(\alpha^2 + \beta^2)  \} \\
  B_{8}  & = & \{ \alpha^4\beta^4,\;
                  \alpha^2\beta^2(\alpha^2 - \beta^2)^2,\;
                  \alpha^3\beta^3(\alpha^2 + \beta^2), \\
       & & \quad  (\alpha^2 - \beta^2)^4,\;
                  (\alpha^2 - \beta^2)^2\alpha\beta(\alpha^2 + \beta^2) \} \\
  & & \dots \\
\end{array}
\end{equation}

\noindent{\bf Example \ref{SectionModesS3D}.2.} For $m = 2$ the
construction yields
\begin{equation}
\begin{array}{lcl}
  B_{0}  & = & \{ 1 \} \\
  B_{2}  & = & \{ \} \\
  B_{4}  & = & \{ \alpha^2\beta^2,\;
                  \alpha^4 + \beta^4 \} \\
  B_{6}  & = & \{ \alpha\beta(\alpha^4 - \beta^4) \} \\
  B_{8}  & = & \{ \alpha^4\beta^4,\;
                  \alpha^2\beta^2(\alpha^4 + \beta^4),\;
                  (\alpha^4 + \beta^4)^2 \} \\
  B_{10} & = & \{ \alpha^3\beta^3(\alpha^4 - \beta^4),\;
                  (\alpha^4 + \beta^4)\alpha\beta(\alpha^4 - \beta^4) \} \\
  B_{12} & = & \{ \alpha^6\beta^6,\;
                  \alpha^4\beta^4(\alpha^4 + \beta^4),\;
                  \alpha^2\beta^2(\alpha^4 + \beta^4)^2,\;
                  (\alpha^4 + \beta^4)^3 \} \\
  & & \dots \\
\end{array}
\end{equation}

\noindent{\bf Example \ref{SectionModesS3D}.3.} For $m = 3$ the
construction yields
\begin{equation}
\begin{array}{lcl}
  B_{0}  & = & \{ 1 \} \\
  B_{2}  & = & \{ \} \\
  B_{4}  & = & \{ \alpha^2\beta^2 \} \\
  B_{6}  & = & \{ \alpha^6 - \beta^6 \} \\
  B_{8}  & = & \{ \alpha^4\beta^4,\;
                  \alpha\beta(\alpha^6 + \beta^6) \} \\
  B_{10} & = & \{ \alpha^2\beta^2(\alpha^6 - \beta^6) \} \\
  B_{12} & = & \{ \alpha^6\beta^6,\;
                  \alpha^3\beta^3(\alpha^6 + \beta^6),\;
                  (\alpha^6 - \beta^6)^2 \}
    \phantom{,\alpha^2\beta^2(\alpha^4 + \beta^4)^2} \\
  & & \dots \\
\end{array}
\end{equation}

For all $k \geq 0$ and all $m > 0$ (including $m = 1$) a
straightforward induction shows the number of base modes $|B_k|$
to be
\begin{equation}
      \left\lfloor\frac{k/2}{2}\right\rfloor
    + \left\lfloor\frac{k/2}{2}\right\rfloor
    + \left\lfloor\frac{k/2}{m}\right\rfloor
    - \frac{k}{2}
    + 1.
\end{equation}
This formula looks superficially different from the result of
Ikeda's dimension counting argument~\cite{Ikeda}, which states
separate formulas for $k \equiv 0$~(mod~4) and $k \equiv
2$~(mod~4), but turns out to be equivalent.  Note, however, that
Ikeda's Theorem~4.3 needlessly excludes the case $m = 1$ and
accidently excludes the eigenmodes for $k = 2m$ when $m$ is odd.

The space $S^3/D_m^*$ admits no modes of odd wavenumber $k$
because such modes fail to be invariant under the antipodal map
$-id \in D_m^*$. \\

\noindent{\bf Theorem \ref{SectionModesS3D}.4.} {\it For each $m >
0$ and each even wavenumber $k$, the polynomials
\begin{equation}
\begin{array}{lll}
    &
    (\alpha^2 \beta^2)^{a}\,
    (\alpha^{2m} + \beta^{2m})^{b}\,
    (\alpha\beta(\alpha^{2m} - \beta^{2m}))^{c}
    & \mathrm{if}\;m\;\mathrm{is~even} \\
    \mathrm{or} & & \\
    &
    (\alpha^2 \beta^2)^{a}\,
    (\alpha^{2m} - \beta^{2m})^{b}\,
    (\alpha\beta(\alpha^{2m} + \beta^{2m}))^{c}
    & \mathrm{if}\;m\;\mathrm{is~odd} \\
\end{array}
\end{equation}
satisfying $a \geq 0$, $b \geq 0$, $c \in \{0,1\}$ and $4a + 2 m b
+ (2m+2)c = k$, together with all their siblings (defined in
Section~\ref{SectionSiblingModes}), form a basis for the
eigenmodes of $S^3/D_m^*$.  The basis has dimension
\begin{equation}
\label{EqnDimensionD}
  \left(
      \left\lfloor\frac{k/2}{2}\right\rfloor
    + \left\lfloor\frac{k/2}{2}\right\rfloor
    + \left\lfloor\frac{k/2}{m}\right\rfloor
    - \frac{k}{2}
    + 1
  \right)
  (k + 1).
\end{equation}
and the polynomials may be easily constructed by the iterative
procedure described above.  There are no modes for odd $k$.}

\vskip0.25cm

\noindent{\it Proof.}  Contained in preceding discussion.
$\blacksquare$ \\

\section{Modes of the Lens Spaces $L(p,1)$}
\label{SectionModesS3Z}

The ``binary cyclic groups'' $Z_n^*$ give only even order lens
spaces, so rather than trying to lift group elements and modes
from $S^2$ to $S^3$ as we did in
Sections~\ref{SectionModesS3I}-\ref{SectionModesS3D}, let us
instead work in $S^3 \subset \mathbb C^2$ directly.  The lens
space $L(p,1)$ is generated by the Clifford translation
\begin{equation}
  \left(
    \begin{array}{c}
      \alpha \\
      \beta \\
    \end{array}
  \right)
  \mapsto
  \left(
    \begin{array}{cc}
      e^{\frac{2\pi i}{p}} & 0 \\
      0 & e^{-\frac{2\pi i}{p}} \\
    \end{array}
  \right)
  \left(
    \begin{array}{c}
      \alpha \\
      \beta \\
    \end{array}
  \right)
  =
  \left(
    \begin{array}{r}
      e^{\frac{2\pi i}{p}} \alpha \\
      e^{-\frac{2\pi i}{p}} \beta \\
    \end{array}
  \right)
\end{equation}
which immediately tells us that the group preserves the mode
$\alpha^{a}\beta^{b}$ if and only if $a - b \equiv 0$~(mod~$p$).
Thus for a given wavenumber $k$ (odd or even) the number of base
modes is simply the number of ways to choose nonnegative $a$ and
$b$ satisfying $a + b = k$ and $a - b \equiv 0$~(mod~$p$).
Examining these constraints on the integer lattice in the $(a,b)$
plane gives the dimension counts stated in
Theorem~\ref{SectionModesS3Z}.1.  These counts agree with
previously obtained real-variable eigenmodes for
lens spaces~\cite{LUW}. \\

\noindent{\bf Theorem \ref{SectionModesS3Z}.1.} {\it For each $p >
0$ and $k \geq 0$ (odd or even), the monomials
\begin{equation}
  \alpha^{a}\beta^{b}
\end{equation}
for $a \geq 0$, $b \geq 0$, $a + b = k$ and $a - b \equiv
0$~(mod~$p$), together with all their siblings (defined in
Section~\ref{SectionSiblingModes}), form a basis for the
eigenmodes of $L(p,1) = S^3/Z_p$ with dimension
\begin{equation}
\begin{array}{cl}
  \left(2\lfloor\frac{k}{ p}\rfloor + 1\right)(k + 1)   & \mathrm{if}~p~\mathrm{is~even~and}~k~\mathrm{is~even} \\
  0                                                     & \mathrm{if}~p~\mathrm{is~even~and}~k~\mathrm{is~odd } \\
  \left(2\lfloor\frac{k}{2p}\rfloor + 1\right)(k + 1)   & \mathrm{if}~p~\mathrm{is~odd ~and}~k~\mathrm{is~even} \\
  2\left(\lfloor\frac{k-p}{2p}\rfloor + 1\right)(k + 1) & \mathrm{if}~p~\mathrm{is~odd ~and}~k~\mathrm{is~odd } \\
\end{array}
\end{equation}
In the special case of the ``binary cyclic groups'' $S^3/Z_n^* =
L(2n,1)$ the dimension equals
\begin{equation}
\label{EqnDimensionZ}
  \left(
      \left\lfloor\frac{k/2}{n}\right\rfloor
    + \left\lfloor\frac{k/2}{n}\right\rfloor
    + 1
  \right)
  (k + 1).
\end{equation}
for even $k$ and zero for odd $k$, in analogy with
Theorems~\ref{SectionModesS3I}.5, \ref{SectionModesS3O}.1,
\ref{SectionModesS3T}.1 and \ref{SectionModesS3D}.4.}

\vskip0.25cm

\noindent{\it Proof.}  Contained in preceding discussion.
$\blacksquare$ \\

\section{Modes of the 3-Sphere $S^3$}
\label{SectionModesS3}

For sake of completeness, let us note that the 3-sphere fits
comfortably within the framework of the preceding sections. \\

\noindent{\bf Theorem \ref{SectionModesS3}.1.} {\it For each $k
\geq 0$ (odd or even), the monomials
\begin{equation}
  \alpha^{a}\beta^{b}
\end{equation}
for $a \geq 0$, $b \geq 0$ and $a + b = k$, together with all
their siblings (defined in Section~\ref{SectionSiblingModes}),
form a basis for the eigenmodes of $S^3$ with dimension $(k +
1)^2$.}

\vskip0.25cm

\noindent{\it Proof.}  The polynomials $\alpha^a \beta^b$ and
their siblings are linearly independent, are harmonic by
Propositions~\ref{SectionHarmonicity}.1 and
\ref{SectionHarmonicity}.2, and meet the well-known dimension $(k
+ 1)^2$ of the mode space of $S^3$.
$\blacksquare$ \\

\section{Future Work}
\label{SectionFutureWork}

The present article provides explicit polynomials for the
eigenmodes of all globally homogeneous spherical 3-manifolds. Four
directions present themselves for future work.

\subsection{Find a unified approach}
\label{SubsectionFutureUnifiedApproach}

Mathematically, the similarity of the base mode dimension
formulas~(\ref{EqnDimensionI}), (\ref{EqnDimensionO}),
(\ref{EqnDimensionT}), (\ref{EqnDimensionD}) and
(\ref{EqnDimensionZ}), along with their obvious dependence on the
underlying $pqr$~triangle group, is too striking to ignore. While
this could be merely an effect of dimension counts, one suspects a
single unified proof might be possible. Here is a plausible
approach, for all cases except the odd order lens spaces. Starting
with one of the base modes, whose twist equals its (even) degree
$k$, apply the $twist_{-}$ operator $k/2$ times to obtain a
sibling of twist $0$. Because the twist is zero, the sibling
projects down from a mode of $S^3$ to a well-defined mode of
$S^2$.  Elementary considerations imply that the mode is pure real
or pure imaginary, so after possibly multiplying by $i$ we may
assume it is pure real.  Even more elementary considerations imply
the result is harmonic on $S^2$. The idea for a unified proof is
to reverse-engineer this process. First construct polynomials on
$S^2$ forming a basis for the modes invariant under a
$pqr$~triangle group, then lift those modes to twist~$0$ modes of
$S^3$, and finally apply both the $twist_{-}$ and the $twist_{+}$
operator $k/2$ times each to obtain all siblings.

\subsection{Orthonormalize the basis}
\label{SubsectionFutureOrthonormalization}

The theorems proved in this paper provide a basis for the space of
eigenmodes of every homogeneous spherical 3-manifold. CMB
simulations, however, require not just an arbitrary basis, but
rather an orthonormal basis.  The author has tried to find a
simple and elegant way to normalize the given bases, but so far
without success.  The problem seems to require a fundamentally new
insight.  As a fallback plan, numerical methods would work, but
with a heavy cost in terms of slower computation and diminished
accuracy.

\subsection{Transfer to spherical coordinates}
\label{SubsectionFutureSphericalCoordinates}

The modes found in this paper are wonderfully simple:  they are
explicit polynomials of fairly low degree in only two variables.
For CMB simulations, however, we need to relate these modes to the
last scattering surface, which in effect requires transferring the
modes to spherical coordinates.  The change-of-coordinates must be
done quickly and accurately;  otherwise the explicit polynomials
will not be significantly faster or better than numerical methods
(e.g. using more than 10000 eigenfunctions \cite{aurich2}).

\subsection{Improve CMB simulations}
\label{SubsectionFutureSimulations}

Cosmologically, the next step is to apply the polynomials we have
found to streamline the cosmological simulations.  The late Jesper
Gundermann began this process and reported improved accuracy, the
elimination of small extraneous imaginary quantities from the
calculation, and the possibility of extending the cosmological
simulations to higher wavenumbers $k$.  A full application of the
polynomials will require satisfactory solutions to the problems
discussed in Sections~\ref{SubsectionFutureOrthonormalization} and
\ref{SubsectionFutureSphericalCoordinates}.

\section*{Acknowledgments}
\label{SectionAcknowledgments}

The author thanks the NSF for its support, and Jesper Gundermann
for suggesting improvements to this article.

\end{document}